\documentclass[12pt,oneside]{amsart}
\usepackage{amssymb}

      \theoremstyle{plain}
      \newtheorem{theorem}{Theorem}[section]

      \theoremstyle{definition}

      \theoremstyle{remark}

 \newcommand{\CC}{{\mathbb C}}
      \newcommand{\RR}{{\mathbb R}}
\newcommand{\ZZ}{{\mathbb Z}}

      \makeatletter
      \def\@setcopyright{}
      \def\serieslogo@{}
      \makeatother

\begin{document}

\date{~}
   \author{Michael Atiyah}
   \address{Trinity College, Cambridge and the University of Edinburgh}
   \email{M.Atiyah@ed.ac.uk}
   \curraddr{~\newline
School of Mathematics\newline
University of Edinburgh\newline
James Clerk Maxwell Building\newline
Peter Guthrie Tait Road\newline
Edinburgh EH9 3FD\newline
Scotland, UK}

   \title{The Non-Existent Complex 6-Sphere}



   \dedicatory{Dedicated to S.S.Chern, Jim Simons and Nigel Hitchin}

   \date{\today}

\begin{abstract}
The possible existence of a complex structure on the 6-sphere has been a famous unsolved problem for over 60 years.  In that time many ``solutions'' have been put forward, in both directions. Mistakes have always been found. In this paper I present a short proof of the non-existence, based on ideas developed, but not fully exploited, over 50 years ago.
\end{abstract}

   \maketitle

   \section{Introduction: history of the problem}

The fact that the 2-sphere $S^2$ can be identified with the complex projective line has been known for centuries. In higher dimensions we have two families of manifolds:
\begin{itemize}
\item[a)] the $2n$-spheres $S^{2n}$
\item[b)] the complex projective spaces $P_n(\CC)$
\end{itemize}

   For $n>1$, these have quite different topologies and so cannot be identified.  It remains possible that the $2n$-sphere has a complex structure quite different from that of $P_n(\CC)$, but if so this cannot have a K\"ahler metric.

  In the early fifties, topology had made big strides and it was possible, using Steenrod squares, to eliminate all values of $n$ except 3, namely the 6-sphere.  This has turned out to be a very hard problem which has resisted determined efforts by the best geometers, most recently S.S.Chern, who in the last year of his life, made some real progress on the problem.  An excellent account of the problem and of Chern's work has been written up by Bryant [13].

  The reason why the $S^6$ problem appears so hard is that $S^6$ has a well-known almost complex structure $J(0)$, which comes from the octonions, when the 6-sphere is viewed as the unit sphere in the 7-space of imaginary octonions.  This however is not integrable, meaning that the operator $\bar{\partial}$ does not have square zero.  This purely algebraic fact comes ultimately from the non-associativity of the octonions.

  Standard topology appears unable to determine whether $S^6$ might have another almost complex structure $J$ which is integrable.  Nothing in the topologist's armoury seems available, notably we have no K\"ahler metric, which is the usual way that complex structures lead to topology via Hodge theory.  It seems clear that some fundamentally new methods are needed, but we will first review the standard theory of the 1950’s which sets the context.

  In 1953, at a conference in Cornell, Hirzebruch [15] listed a large number of important problems in Geometry.  Many of these were solved in the subsequent decade, notably those centring round the Hirzebruch-Riemann-Roch Theorem HRR.  A major step was taken by Atiyah-Singer with the index theorem [8], which shifted the arena from projective algebraic geometry to differential geometry.  In particular, for complex analytic manifolds, the restriction to K\"ahler manifolds was removed.  The significance of this was immediately recognised by Kodaira who was then able to complete the coarse classification of compact complex surfaces, including such non-K\"ahler surfaces as $S^1 \times  S^3$.

  Shortly after HRR a major breakthrough in algebraic geometry was made by Grothendieck [14] with the introduction of $K$-theory.  This inspired Atiyah-Hirzebruch [6] to develop a topological analogue based on the Bott periodicity theorem for the unitary group [12], and this K-theory replaced cohomology as the natural home of index theory.  Now Bott periodicity had counterparts for the orthogonal and symplectic groups, with period 8 and ``semi-period''  4, switching orthogonal and symplectic.  It also applied to spin-manifolds. This naturally tied up with index theory over the reals [10] where the fundamental operator is the Dirac operator and the mod 2 invariants appeared as the mod 2 dimensions of the null spaces of skew-adjoint operators [9].

  Indices which are integers can be computed via Chern classes as rational cohomology or integrals.  Mod 2 invariants are much more elusive, except in the lowest dimensions, where they come from the first and second Stiefel-Whitney classes.  However, with $KO$-theory, they can be handled effectively.

  Let us first re-examine the success of index theory in giving a proof of HRR for non-K\"ahler complex manifolds.  It hinges on the following facts:

1.1 The holomorphic Euler characteristic (the arithmetic genus) uses the integer grading of sheaf cohomology, while the index uses only the parity of $K$-theory.

1.2 When we replace the Dolbeault complex by the single elliptic operator $D$ given by adding $\bar{\partial}$ to its adjoint, with respect to an arbitrary positive Hermitian metric, the fact that $DD^*$ does not preserve the grading is irrelevant.  It can be deformed away without affecting the ellipticity.

1.3 Thus the index only uses the almost complex structure and we get the integrality of the Todd genus.

1.4 Similarly, using a spin structure in the right dimensions, captures the $KO$ mod 2 index which is defined as the mod 2 dimension of the null space of the Dirac operator.
For example, in dimension 2, the $KO$ (spin) mod 2 index goes back to Riemann and was studied in [2].

1.5 Another example in dimension 6 arises from a rank 2 bundle on the complex projective 3-space. A mod 2 invariant in this situation was studied in [7].

    Recall finally that Atiyah-Singer theory does not use any special metric, it uses an arbitrary auxiliary Riemannian metric leading to elliptic operators. Its main theorem is then independent of the metric and asserts that analytical indices are equal to topological indices.  This applies to both integer and mod 2 invariants.
So there is no need to distinguish between analysis and topology in this whole area. This provides a new tool whose power so far has only been fully exploited by Kodaira.

It is this tool, in the mod 2 context, which we will exploit to tackle the problem of $S^6$. The version of $K$-theory needed is $KR$ theory, which was developed in detail in [1], and whose key features will be recalled in the next section.  It was motivated both by real algebraic geometry and by real differential operators.

\section{$KR$ Theory}

$KR$ theory is defined for manifolds with an involution, called conjugation, as in algebraic geometry. It has a double index notation $(p,q)$, based on Clifford algebras of quadratic forms of signature $(p,q)$,where the conjugation is changing the sign of the set of $q$ variables. Suspension is given by taking the $KR$-theory with compact supports of the product with $\RR(p,q)$. The $q$ variables give the standard suspension while the p variables give a twisted suspension, in that conjugation multiplies by $-1$ in the fibre.

The complex field $\CC$ has its standard conjugation so that $\CC = \RR(1,1)$.
There is a natural generator $b$ of $KR (1,1)$ and the Bott Periodicity Theorem now asserts that the tensor product with $b$ gives an isomorphism
                   $$ KR (p,q)  \cong  KR(p+1,q+1)$$
This holds over any space with involution, generalizes to the equivariant context, and is a powerful theorem despite its apparent simplicity.  It shows that, up to canonical isomorphism, $KR(p,q)$ depends only on $s=p - q$ and that it is periodic in $s$ with period 8, and semi-periodic with period 4 if we switch orthogonal and symplectic.

Applied when the base space is a point and $(p,q) = (7,1)$ we get an isomorphism
$$KSp(\RR^6) \to KR^{(7,1)}(\hbox{point})  = \ZZ/2$$

    The LHS is the stable homotopy group determined by Bott and the RHS is given by the algebra of Clifford modules. The homomorphism is the natural one as described in [3].

The left hand side is just the reduced $KSp$ of $S^6$ and elements of it are either even or odd (see section 3), and their values under the isomorphism above are respectively 0 or 1.

There are natural forgetful maps from complex $K$-theory to $KR$ theory and in dimension 6 the integers go to 0, so a hypothetical complex structure on $S^6$ would give an even element. But we know one almost complex structure $J(0)$ which is odd, and
{\bf crucially we know that being even or odd is a topological property} and so shared by all almost complex structures (see section 3).   This is a contradiction and so proves

\begin{theorem}
     There is no complex structure on the 6-dimensional sphere
\end{theorem}

\noindent Comments
\smallskip

\noindent 1. The non-trivial part of the proof (in bold lettering in the text above) rests on the Atiyah-Singer index theorem, in parallel to the way it was used by Kodaira

\noindent 2. Since our invariant was a mod 2 index we had to use $KR$ theory (designed for such purposes)

\noindent 3. Hidden in the technical details of $KR$ theory, the mod 2 invariant is really a mod 4 invariant, with non-trivial ``odd'' elements being strictly of order 4 (and  not of order 2). An algebraic geometer, working over the real field knows that one should not distinguish between $i$ and $-i$ whereas one can distinguish between 1 and $-1$. Because of the Atiyah-Singer theory this applies in topology as well (a stronger form of Serre's GAGA).

\noindent 4. For the non-K\"ahler complex surfaces that Kodaira needed, the natural metric was not K\"ahler but had signature (2,2). In our case the natural 8-manifold is not the octonion space of signature (8,0) but its Minkowski counterpart with signature (7,1). The 6-sphere then appears naturally as the base of the light-cone. This can be compared with the standard Minkowski space of signature (3,1), where the 2-sphere is famously the base of the light cone, but there is one important difference. The quaternions are  associative while the octonions are not (this is purely algebraic and independent of signatures),   explaining why the natural almost complex structure on the 2-sphere is integrable, while it is not for the 6-sphere.

\noindent 5. Note that the theorem is a global one.  Locally,   on the punctured sphere, it is false. The theorem has to come from some global ``cohomological'' obstruction. In fact, as we have seen, the cohomology is $KR$ theory.

There are a number of subtleties about almost complex structures which lurk in the background of note 5 and in the details of $KR$ theory.  These deserve separate treatment in the following section.

\section{Almost Complex Structures}

    Bryant \cite{[13]} has given a very clear account of the linear algebra and
the associated differential geometry of almost complex manifolds, which
is precisely what we need. The only difference is that, whereas Bryant
works in a Riemannian context, we will work in a Lorentzian one.  As
noted in comment 2 of section 2, the 6-sphere is the base of the light
cone in Minkowski space  R(7,1) and has a natural conformal structure
but no preferred metric.

   The reason this does not cause problems is because the Atiyah-Singer
theory (which relates analysis to topology) does not need a preferred
Riemannian metric.  This was precisely why Kodaira was able to handle
complex non-K\"ahler surfaces.

   On $S^6$ the Octonions of signature (7,1) induce an almost complex
structure.  In fact, as noted in comment 3 of section 2, we get a
conjugate pair of almost complex structures.  This should be called a
real ACS. Such a conjugate pair can be distinguished locally on a
punctured sphere, but this distinction may not persist globally.  There
are just 2 topological possibilities: either the distinction is global
(called the even case), with the distinction carrying over the puncture,
or it is not (called the odd case).

   The distinction between odd and even elements of $KSp({\mathbb R}^6)$ noted in
section 2 does agree with the distinction made above by virtue of the
Atiyah-Singer index theory, which was essentially designed for the purpose

    For real quadrics the real ACS depends in an interesting way on the
signature of the quadratic form.  Signature (7,1) leads to an odd real
ACS while signatures (5,3) and (3,5) lead to the even case of genuine
complex structures on the complex projective 3-space $P$ or its dual $P^*$.
This all fits naturally into the twistor theory of $S^6$ where triality
tries to seduce us into believing in a mythical complex structure on
$S^6$, but it is well-known that one of the three groups involved in
triality has to be lifted to spin. This is clear for the  algebra of
octonions but, because the obstruction is  captured by a topological
obstruction in a $KR$-group, it cannot be avoided by any deformation.

  Our Theorem can now be interpreted as saying that, on the 6-sphere,
any real ACS is of odd type.  Hence there is no real ACS of of even
type.  An integrable complex structure would define a real ACS of even
type and so cannot exist. The integrability condition is essentially replaced by an  equivalent topological condition.  This for the 6-sphere is precisely what we expected, since there is no local obstruction and we needed a global cohomological obstruction in an appropriate theory.  That theory is just $KR$ theory.

                       A short history of the 6-sphere problem follows as section 4.

\section{History of the problem}

\noindent Ehresmann 1947: Introduced the notion of almost complex
structure and showed that the 6-sphere admits an almost
complex structure, but explicitly points out that he does not
know whether it has a complex structure.\\
Hopf 1947: Proved that $S^4$ and $S^8$ do not admit almost
complex structures.\\
Kirchhoff 1947: Uses octonions to construct an explicit
almost complex structure on $S^6$.\\
 Eckmann-Frohlicher and Ehresmann-Liberman 1951:
Independently prove that Kirchhoff's almost complex
structure on $S^6$ is not integrable to a complex structure.\\
 Borel and Serre 1953: Prove that $S^{2n}$ admits an almost
complex structure if and only if $n = 1$ or 3.\\
Hirzebruch 1954 and Liberman 1955: Remarks that it is still
not known whether $S^6$ has a complex structure.

\section*{Acknowledgments}

   I gratefully acknowledge the help, advice and criticism of many colleagues.  First there are the three ``musketeers'' of the gang of 4, Bott, Hirzebruch and Singer.  Then there is Jim Simons, for whom this paper was written in celebration of his honorary degree from the University of Edinburgh.  As a bonus the first draft was completed on 4 July 2016. Dennis Sullivan and Claude Lebrun showed some healthy and fruitful scepticism.  My younger colleagues, Nigel Hitchin and J\"urgen Berndt, provided expertise and constructive criticism.  Robert Bryant who knows the problem well [13] pushed me hard and, supported by Blaine Lawson, forced me to clarify my argument. Finally my amanuensis Andrew Ranicki has always been at my elbow, speeding me along.

\end{document}